\documentclass[a4paper,reqno,12pt]{amsart}
\usepackage{amsmath,amscd,amssymb,upgreek} 
\usepackage[all]{xy}
\usepackage{graphicx}
\usepackage{ifthen}
\usepackage{yhmath}
\usepackage{enumerate}%\parindent3mm 
\usepackage{pdfsync}
\usepackage{color}
\usepackage{hyperref}

\usepackage[applemac] {inputenc}

\newcommand{\application}[5]{
\ifthenelse{\equal{#1}{0}}{}{#1:}\begin{array}[t]{ccl}
	#2 & \longrightarrow & #3
	\ifthenelse{\equal{#4}{0}}{}{ \\ #4 & \longmapsto & #5} 
\end{array}}

\newcommand{\Theoreme}[1]{\begin{The}#1\end{The}}
\newcommand{\Proposition}[1]{\begin{Prop}#1\end{Prop}}
\newcommand{\Lemme}[1]{\begin{Lem}#1\end{Lem}}
\newcommand{\Corollaire}[1]{\begin{Cor}#1\end{Cor}}
\newcommand{\Definition}[1]{\begin{Def}#1\end{Def}}
\newcommand{\Preuve}[1]{\noindent \textbf{Proof:} #1~\hfill$\blacksquare$\\}

\def\R{\mathbb{R}}
\def\C{\mathbb{C}}
\def\H{\mathbb{H}}

\def\Q{\mathbb{Q}}

\def\Z{\mathbb{Z}}

\def\S{\mathbb{S}}

\begin{document}
\newtheorem{The}{Theorem}[section]
\newtheorem{Prop}[The]{Proposition}
\newtheorem{Def}[The]{Definition}
\newtheorem{Lem}[The]{Lemma}
\newtheorem{Cor}[The]{Corollary}
\newtheorem{Rem}[The]{Remark}
\newtheorem{Exple}[The]{Example}
\title[Line bundles and algebraic dimension of OT-manifolds]{Holomorphic line bundles over domains in Cousin groups and the algebraic dimension of OT-manifolds \footnote{AMS classification 2010: 32J18, 32Q57, 32J99}}
\begin{author}
{L. Battisti \& K. Oeljeklaus}
\end{author}
\begin{abstract}
In this paper we extend results due to Vogt on line bundles over Cousin groups to the case of domains stable by the maximal compact subgroup. This is used in the sequel to show that the algebraic dimension of OT-manifolds is zero. In the last part we establish that certain Cousin groups, in particular those arising from the construction of OT-manifolds, have finite-dimensional irregularity.
\end{abstract}
\begin{address}
{\parbox{10cm}{Laurent BATTISTI \& Karl OELJEKLAUS\\ 
AIX-MARSEILLE UNIVERSITÉ\\  CNRS - %\\ %Centrale Marseille\\ 
LATP - UMR 7353\\F-13453 MARSEILLE \\France\\
\textbf{\url{battisti@cmi.univ-mrs.fr}, \url{karloelj@cmi.univ-mrs.fr}\\}}}
\end{address}
\maketitle
%\tableofcontents
%\nocite{Huckleberry:1983aa}\nocite{Cousin:1910aa}
\vspace{-0.5cm}\section{Introduction}
\noindent A connected complex Lie group which admits no non-constant holomorphic functions is called a \emph{Cousin group}, or also a \emph{toroidal group} in the older literature. These groups are named after P.~Cousin (see \cite{Cousin:1910aa}). If $C$ is an $n$-dimensional Cousin group, it is abelian and can be realized as the quotient of $\C^{n}$ by a lattice $\Lambda$ of $\C^{n}$ of rank $n+m$, with $1\leq m \leq n$.~\\
\indent It is a fact that the line bundle associated to a positive divisor on a compact K\"ahler manifold
has a non-trivial real Chern class. In \cite{Vogt:1983aa}, C. Vogt proved the following in the non-compact context: Let $L$ be a topologically trivial line bundle on a Cousin group. Then $L$ is holomorphically trivial if and only if it admits a non-trivial holomorphic section. 
This result had already been proved by Cousin \cite{Cousin:1910aa}, in the special case $m=1$. His original proof was rewritten by A.T.~Huckleberry and G.A.~Margulis in \cite{Huckleberry:1983aa}, who used it to establish a theorem on
hypersurfaces in quotients of semi-simple complex Lie groups, see also \cite{Berteloot:1988fk}, \cite{Gilligan:2011uq}.\\
% He also established that every line bundle on a Cousin group $C$ comes from a theta factor if and only if the space $H^{1}(C,\mathcal{O})$ is finite-dimensional. ~\\
\indent In the first section of the present paper we generalize Vogt's result to the case of a \emph{domain} $U$ in a Cousin group $C=\C^{n}/\Lambda$ whose inverse image in $\C^{n}$ is convex. This assumption implies in particular that $U$ is invariant under the action of the maximal compact subgroup of $C$.~\\
\indent In the next section we use the previous result to confirm the non-existence of
complex hypersurfaces in OT-manifolds.
In particular, the algebraic dimension of these manifolds is always zero. This result was proved in a special case by L.~Ornea and M.~Verbitsky in \cite{Ornea:2011fk}. OT-manifolds were introduced in \cite{Oeljeklaus:2005aa} and we will recall their construction later.\\
\indent The last section is devoted to
the observation that if the coordinates of the
lattice points of a Cousin $C$ group are all algebraic
integers then the irregularity $\dim_{\C} H^{1}
(  C, \mathcal O)$ is finite.\\

This article is organized as follows. In the first section we recall the notations of \cite{Vogt:1982aa} and \cite{Vogt:1983aa}, and prove the following:~\\~\\
\noindent\textbf{Theorem~\!\ref{propVogtGeneralisee}.} \emph{Let $U\!$ be an open subset of a Cousin group $C\!\cong\!\C^{n}/\Lambda$ whose inverse image $\widetilde{U}$ in $\C^{n}$ is a convex domain. Let $L$ be a topologically trivial holomorphic line bundle over $U$. One has $H^{0}(U,L) \neq 0$ if and only if $L$ is holomorphically trivial.}~\\

In the second section, we first recall the construction of OT-manifolds and then prove:~\\

\noindent\textbf{Theorem~\ref{dimAlgCasGeneral}.} \emph{Let $X$ be an OT-manifold. Then there are no complex-analytic hypersurfaces on $X$, in particular the algebraic dimension of $X$ is zero.}~\\

Before stating the result of the last section we recall:~\\~\\
\noindent\textbf{Theorem~\ref{thCaractClasseCousinVogt} (\cite{Vogt:1983aa}).} \emph{Let $C=\C^{n}/\Lambda$ be a Cousin group. Then the following conditions are equivalent:
\begin{itemize}
\item[1.] The space $H^{1}(C,\mathcal{O})$ is finite-dimensional.
\item[2.] Let $P=(I_{n}~S)$ be a period basis of $\Lambda$. Then there exist constants $C>0$ and $a\geqslant 0$ such that $\|^{t}\sigma S+\,^{t}\tau\|\geqslant C \exp(-a|\sigma|)$ for all $\sigma \in \Z^{n}\setminus\{0\}$ and all $\tau\in\Z^{m}$, where $n+m$ is the rank of $\Lambda$.
\item[3.] Every line bundle over $C$ comes from a theta factor.
\end{itemize}}
~\\
\indent Applying a generalization of Liouville's theorem we obtain 
the following which applies in particular to Cousin groups arising in the construction of OT-manifolds: ~\\

\noindent\textbf{Theorem~\ref{grpeCousinAlgebriqueCondition2}.} \emph{Let $\Lambda \subset \C^{n}$ be a lattice such that $C=\C^{n}/\Lambda$ is a Cousin group with a period basis whose coefficients are all algebraic numbers, then $C$ satisfies the equivalent conditions of theorem~\ref{thCaractClasseCousinVogt}.}

 \section{Preliminaries}

Let us consider a domain $U$ of a Cousin group $C \cong \C^{n}/\Lambda$ whose inverse image $\widetilde{U}=\pi^{-1}(U)$ in $\C^{n}$ is a convex domain, where $\pi: \C^{n} \rightarrow \C^{n}/\Lambda$ is the quotient map. In particular, $U$ is invariant under the action of $\R\Lambda/\Lambda$, admits no non-constant holomorphic functions and $\widetilde{U}$ is invariant under the action of $\R\Lambda$ and is evidently Stein. We recall notations, definitions and results from \cite{Vogt:1982aa} which adapt directly to our situation.

\Definition{A map  $\alpha: \Lambda \times \widetilde{U} \rightarrow \C^{*}$ is called a \textbf{factor of automorphy} if it satisfies the following properties:
\begin{enumerate}
\item[a)] $\alpha_{\lambda}: \widetilde{U} \rightarrow \C^{*}$,  $\alpha_{\lambda}(z):= \alpha(\lambda,z)$ is holomorphic for all $\lambda \in \Lambda$,
\item[b)] $\alpha(0,z)=1$ for all $z \in \widetilde{U}$,
\item[c)] $\alpha(\lambda+\lambda',z)=\alpha(\lambda,z+\lambda')\alpha(\lambda',z)$ for all $\lambda, \lambda' \in \Lambda$ and all $z\in \widetilde{U}$.
\end{enumerate}}

If $\alpha$ is a factor of automorphy, for every $\lambda \in\Lambda$ there exists a holomorphic function $a_{\lambda}: \widetilde{U}\rightarrow \C$ unique up to an additive constant $2i\pi k_{\lambda}$ with $k_{\lambda}\in\Z$, such that $\alpha(\lambda,z)=\exp(a_{\lambda}(z))$.  

\Definition{\label{conditionCDefSommant}A map $a: \Lambda \times \widetilde{U} \rightarrow \C$ is called a \textbf{summand of automorphy} if the following three conditions are satisfied: 
\begin{enumerate}
\item[a)] $a(\lambda, ~\cdot~): \widetilde{U}\rightarrow \C$ is holomorphic for all $\lambda \in\Lambda$, %where we set $a_{\lambda}(z):=a(\lambda,z)$,
\item[b)] $a(0,z) = 0$ for all $z\in \widetilde{U}$,
\item[c)] $a(\lambda+\lambda',z)=a(\lambda,z+\lambda')+a(\lambda',z)$ for all $\lambda, \lambda' \in \Lambda$ and all $z\in\widetilde{U}$.
\end{enumerate}}

%If $L \overset{p}{\rightarrow} U$ is a holomorphic line bundle, then the pull-back 
%$$\pi^{*}L:=\widetilde{U}\times_{U} L = \{(z,v) \in \widetilde{U}\times L ~|~\pi(z) = p(v)\}$$
%is a holomorphic line bundle over $\widetilde{U}$, therefore analytically trivial because $\widetilde{U}$ is Stein. Let $\application{\varphi}{\pi^{*}L}{\widetilde{U}\times \C}{(z,v)}{(z,\varphi_{z}(v))}$ be a trivialization of $\pi^{*}L$. Then the map $\alpha$ defined by $\application{\alpha}{\Lambda\times \widetilde{U}}{\C^{*}}{(\lambda,z)}{\varphi_{z+\lambda}\varphi_{z}^{-1}}$ is a factor of automorphy.
%
%Si $\application{\psi}{\pi^{*}L}{\widetilde{U}\times \C}{(z,v)}{(z,\psi_{z}(v))}$ est une autre trivialisation de $\pi^{*}L$, l'application $h: \widetilde{U}\rightarrow \C^{*}$ qui à $z$ associe $\psi_{z}\varphi_{z}^{-1}$ est holomorphe. Si $\beta$ est le factor of automorphy associé à $\psi$, alors on a l'égalité 
%\begin{equation}\label{equivFactAutomorph}
%\beta(\lambda,z) = h(z+\lambda)\alpha(\lambda,z)h^{-1}(z).
%\end{equation}
%Ceci nous amène à introduire la définition suivante:

\Definition{\label{DefEquivSommantsEtFacteurs}Two factors of automorphy $\alpha$ and $\beta$ are called \textbf{equivalent} if there is a holomorphic function $h: \widetilde{U}\rightarrow \C^{*}$ satisfying 
\begin{equation*}\label{equivFactAutomorph}
\beta(\lambda,z) = h(z+\lambda)\alpha(\lambda,z)h^{-1}(z).
\end{equation*}
Two summands of automorphy are called \textbf{equivalent} if the induced factors of automorphy are equivalent.}

%À chaque fibré en droites $L$ sur $U$ on a associé une unique classe d'équivalence de facteurs d'automorphie. Si $\alpha: \Lambda \times \widetilde{U}\rightarrow \C^{*}$ est un factor of automorphy, on construit un fibré en droites $L$ sur $U$ en regardant le quotient de $\widetilde{U}\times \C$ par l'action de $\Lambda$ définit comme suit:
%
%$$\lambda.(z,v):= (z+\lambda, \alpha(\lambda,z)v)$$
%pour $\lambda \in \Lambda, z \in \widetilde{U}$ et $v\in\C$.

Since $\widetilde{U}$ is Stein, one has the following proposition (see \cite{Vogt:1982aa}):
\Proposition{\label{bijSommantsClassesEquiv}There is a bijection between equivalence classes of factors of automorphy and equivalence classes of line bundles over $U$.}

%Indeed, if $\alpha: \Lambda \times \widetilde{U}\rightarrow \C^{*}$ is a factor of automorphy, one gets a holomorphic line bundle $L$ over $U$ by taking the quotient of $\widetilde{U}\times \C$ under the following action of $\Lambda$:
%
%$$\lambda.(z,v):= (z+\lambda, \alpha(\lambda,z)v)$$
%for $\lambda \in \Lambda, z \in \widetilde{U}$ and $v\in\C$.
%
%Conversely, if $L \overset{p}{\rightarrow} U$ is a holomorphic line bundle, then the pull-back 
%$$\pi^{*}L:=\widetilde{U}\times_{U} L = \{(z,v) \in \widetilde{U}\times L ~|~\pi(z) = p(v)\}$$
%is a holomorphic line bundle over $\widetilde{U}$, therefore analytically trivial because $\widetilde{U}$ is Stein and convex. Let $\application{\varphi}{\pi^{*}L}{\widetilde{U}\times \C}{(z,v)}{(z,\varphi_{z}(v))}$ be a trivialization of $\pi^{*}L$. Then the map $\alpha$ defined by $\application{\alpha}{\Lambda\times \widetilde{U}}{\C^{*}}{(\lambda,z)}{\varphi_{z+\lambda}\varphi_{z}^{-1}}$ is a factor of automorphy.\\

We recall now two normal forms for the lattice $\Lambda$ defining a Cousin group  (see \cite{Vogt:1982aa}, propositions~1 and 2): 
\Proposition{\label{formesNormalesReseauCousin}Let $C=\C^{n}/\Lambda$ be a Cousin group. Then the real rank of $\Lambda$ is $r=n+m$ with $1\leqslant m \leqslant n$ and:
\begin{itemize}
\item[1.] The Cousin group $C$ has a period basis, i.e. a basis of $\Lambda$, of the form $P=(I_{n}~S)$; a lattice $\Lambda$ defined by such a matrix gives a Cousin group if and only if $^{t}\sigma S \not\in\Z^{m}$ for all $\sigma \in \Z^{n}\setminus\{0\}$.
\item[2.] The Cousin group $C$ has a period basis of the form $P=\left(\begin{array}{lr} 0& T\\I_{n-m}& R
\end{array}\right)$ where $T=(I_{m}~S)$ is the period basis of a complex torus of dimension $m$ and $R$ is a real matrix; a lattice $\Lambda$ defined by such a matrix gives a Cousin group if and only if $^{t}\sigma R \not\in\Z^{2m}$ for all $\sigma \in \Z^{n-m}\setminus\{0\}$.
 \end{itemize}}

\subsection{Vogt's results generalized}Given a factor of automorphy on $\widetilde{U}$, one gets a line bundle $L$ on $U$ by %the remark following 
proposition~\ref{bijSommantsClassesEquiv}. The following proposition characterizes the topologically trivial bundles obtained that way. Its proof is the same as in the case of a Cousin group, see for instance \cite{Abe:2001aa}, Lemma~2.1.9.
\Proposition{\label{sommantEtTopTrivial}Let $L \overset{p}{\rightarrow} U$ be a holomorphic line bundle given by a factor of automorphy $\alpha: \Lambda \times \widetilde{U} \rightarrow \C^{*}$. % and for all $\lambda \in \Lambda$, let
%\[\application{a_{\lambda}}{\widetilde{U}}{\C}{z}{a(\lambda, z)}\] be a holomorphic function such that $\alpha(\lambda,z)= \exp(a(\lambda,z))$. 
Then the following two conditions are equivalent: 
\begin{enumerate}
\item[a)] the bundle $L$ is topologically trivial,
\item[b)] there is a summand of automorphy $a: \Lambda \times \widetilde{U} \rightarrow \C$ such that $\alpha(\lambda,z) = \exp(a(\lambda,z))$.
%\item[c)] $a(\lambda,z+\lambda')+a(\lambda',z)=a(\lambda',z+\lambda)+a(\lambda,z)$ for all $\lambda, \lambda' \in \Lambda$ and $z \in \widetilde{U}$.
\end{enumerate}
}

The previous proposition now gives a first normal form in the class of a summand of automorphy:

\Proposition{Let $\Lambda \subset \C^{n}$ be a lattice whose period basis is given in the form $(I_{n}~S)$. Then every summand of automorphy $b: \Lambda\times \widetilde{U} \rightarrow \C$ is equivalent to a summand of automorphy $a: \Lambda \times \widetilde{U}\rightarrow\C$ with the following properties: 
\begin{enumerate}
\item[a)] $a(\sigma,z) = 0$ for all $\sigma \in \Z^{n}$ and all $z\in\widetilde{U}$,
\item[b)] for all $\lambda \in \Lambda$, the holomorphic function $\displaystyle \application{a_{\lambda}}{\widetilde{U}}{\C}{z}{a(\lambda,z)}$ is $\Z^{n}$-periodic.
\end{enumerate}}
\Preuve{The projection $\pi: \C^{n} \rightarrow C$ factors through $\C^{n}/\Z^{n}$: \[\C^{n} \overset{\pi_{2}}{\longrightarrow} \C^{n}/\Z^{n} \overset{\pi_{1}}{\longrightarrow} C,\]
with $\pi = \pi_{1}\circ\pi_{2}$. The image of $\widetilde{U}$ by $\pi_{2}$ in $\C^{n}/\Z^{n}\cong (\C^{*})^{n}$ is Stein by Grauert-Docquier's theorem. Let $L$ be the line bundle defined by the summand of automorphy $b$. By the previous proposition, $L$ is topologically trivial, so is $\pi_{1}^{*}L$. Because every topologically trivial line bundle over a Stein manifold is holomorphically trivial, one obtains a trivialization of $\pi^{*}L$ by taking the pull-back of a trivialization of $\pi_{1}^{*}L$. The summand of automorphy $a: \Lambda \times \widetilde{U}\rightarrow \C$ given by this trivialization is equivalent to $b$ in the sense of definition~\ref{DefEquivSommantsEtFacteurs} and it satisfies condition $a)$. Condition $b)$ is a consequence of both $a)$ and the definition of a summand of automorphy.}
%\Remarque{The previous proof is a direct adaptation of Vogt's proof; the only modification is that one has to ensure that $\pi_{2}(\widetilde{U})$ remains Stein. }

One also has a second normal form for summands of automorphy:
\Proposition{\label{formNorm2e}Let $\Lambda \subset \C^{n}$ be a lattice whose basis is given under the form \[P=\left(\begin{array}{lr} 0& T\\I_{n-m}& R
\end{array}
\right),\] such that $C=\C^{n}/\Lambda$ is a Cousin group. Then every summand of automorphy $b: \Lambda \times \widetilde{U} \rightarrow \C$ is equivalent to a summand of automorphy $a: \Lambda\times \widetilde{U}\rightarrow \C$ with the following properties: 
\begin{enumerate}
\item[a)] $a(\lambda,z) = a(\lambda,z_{m+1}, ..., z_{n})$ for all $\lambda \in \Lambda$ and all $z\in\widetilde{U}$,
\item[b)] $a(\lambda,z) = 0$ for all $\lambda \in \left(\begin{array}{c}0\\\Z^{n-m}\end{array}
\right)$ and all $z\in\widetilde{U}$,
\item[c)] for all $\lambda \in \Lambda$, the holomorphic function $\displaystyle \application{a_{\lambda}}{\widetilde{U}}{\C}{z}{a(\lambda,z)}$ is $\left(\begin{array}{c}0\\\Z^{n-m}\end{array}
\right)$-periodic.
\end{enumerate}}
The proof given in \cite{Vogt:1982aa} readily adapts to our case so we will not repeat it here.

Now, we state the main result of this section:
\Theoreme{\label{propVogtGeneralisee}Let $U$ be an open subset of a Cousin group $C \cong \C^{n}/\Lambda$ whose inverse image $\widetilde{U}$ in $\C^{n}$ is a convex domain. Let $L$ be a topologically trivial holomorphic line bundle over $U$. One has $H^{0}(U,L) \neq 0$ if and only if $L$ is holomorphically trivial.}
\Preuve{The sufficient condition is clear so we prove the necessary condition. We choose the second normal form of the lattice of proposition~\ref{formesNormalesReseauCousin}. Since $L$ is topologically trivial, proposition~\ref{sommantEtTopTrivial} implies that it is given by a summand of automorphy $a:\Lambda \times \widetilde{U} \rightarrow \C$. 
By using the normal form of proposition~\ref{formNorm2e} we may assume that:% the two following conditions hold:
\begin{itemize}
\item[$\boldsymbol{\cdot}$]$a(\lambda,z) = a(\lambda,z_{m+1}, ..., z_{n})$ for all $\lambda \in \Lambda$ and 
\item[$\boldsymbol{\cdot}$]$a_{\lambda}: \widetilde{U}\rightarrow \C$, $z\mapsto a(\lambda,z)$ is $\left(\begin{array}{c}0\\\Z^{n-m}\end{array}
\right)$--periodic for all $\lambda$.
\end{itemize}

Now consider the Fourier series of  $a_{\lambda}$:  
\[a_{\lambda}(z) = \!\!\!\!\sum_{\sigma\in\Z^{n-m}}\!\!\!\!a_{\lambda,\sigma}e^{2i\pi ^{t}\sigma(z_{m+1},...,z_{n})}.\]
%where $\langle ~\cdot~,~\cdot~ \rangle$ denotes the usual scalar product on the Euclidian space.
Since $a$ is a summand of automorphy, we have: 
\[a_{\lambda}(z+\lambda') + a_{\lambda'}(z) =a_{\lambda'}(z+\lambda) + a_{\lambda}(z) \]
hence: \[a_{\lambda,\sigma}(1-e^{2i\pi ^{t} \sigma \lambda_{2}' })=a_{\lambda',\sigma}(1-e^{2i\pi ^{t} \sigma \lambda_{2}})\] for $\lambda, \lambda' \in \Lambda$ (where $\lambda_{2}$ [resp. $\lambda_{2}'$] is the $(n-m)$-tuple consisting of the $n-m$ last coordinates of $\lambda$ [resp. $\lambda'$]) and $z\in\widetilde{U}$.

We want to extend the map $a$ to $\R \Lambda \times \widetilde{U}$ such that
\begin{equation}\label{eqvoulue}
a(x+x',z)=a(x,z+x')+a(x',z) \text{  for all } x,x'\in \R\Lambda, z\in\C^{n}.
\end{equation}

First, we define $a(x, z)$ for an element $x$ of $\R\Lambda$ of the form $x = r \lambda$ by:
\[a(x,z):= \!\!\!\!\sum_{\sigma \in \Z^{n-m}}\!\!\!\! a_{\lambda,\sigma}C(\sigma,r,\lambda)e^{2i\pi ^{t} \sigma(z_{m+1},...,z_{n})},\]

where \[C(\sigma,r,\lambda) = \left\{
\begin{array}{cc}
\displaystyle \frac{1-e^{2i\pi r ^{t} \sigma \lambda_{2}}}{1-e^{2i\pi ^{t} \sigma \lambda_{2} }} & \text{if }^{t}\sigma \lambda_{2}\not\in\Z\phantom{.} \\
&\\
\displaystyle \lim_{u\rightarrow0} \frac{1-e^{iru}}{1-e^{iu}}=r & \text{if }^{t}\sigma \lambda_{2}\in\Z.
 \end{array}\right.\]
%If $\langle \sigma,\lambda_{2} \rangle\in\Z$, then the quotient in this definition is replaced by the quantity $$ \lim_{t\rightarrow0} \frac{1-e^{irt}}{1-e^{it}} = r.$$

The Fourier series defining $a(x,z)$ is convergent since for any fixed real $r$, the factors $C(\sigma,r,\lambda)$ are uniformly bounded by $r$.

Now for an element of $\R\Lambda$ of the form $x = \displaystyle \sum_{j=1}^{n+m}r_{j}\lambda_{j}$ we define the differentiable function $a(x,z)$ everywhere on $\R\Lambda \times \widetilde{U}$ using equation (\ref{eqvoulue}).

Consider a section $s$ of $L$. We shall prove that if $s$ vanishes at one point, then it is constant equal to zero. This section corresponds to a holomorphic function $f: \widetilde{U}\rightarrow \C$, satisfying 
\begin{equation}\label{condSurF}
f(z+\lambda) = e^{a(\lambda,z)}f(z)
\end{equation}
 for all $\lambda \in \Lambda$ and all $z\in\widetilde{U}$.

Let $z_{0} \in \widetilde{U}$, consider the function: \[\application{t_{z_{0}}}{\R\Lambda}{\C}{x}{f(z_{0}+x)e^{-a(x,z_{0})}.}\] 
Then, for all $\lambda \in \Lambda$ and all $x\in\R\Lambda$, one has:
\begin{eqnarray*}
t_{z_{0}}(x+\lambda) &=& f(z_{0}+x+\lambda)e^{-a(x+\lambda,z_{0})}\\ 
&=& f(z_{0}+x)e^{a(\lambda,z_{0}+x)-a(x+\lambda,z_{0})}\\
&=& f(z_{0}+x)e^{-a(x,z_{0})}\\
&=&t_{z_{0}}(x).
\end{eqnarray*}
Hence $t_{z_{0}}$ is a $\Lambda$-periodic function and therefore bounded by a constant $C_{z_{0}}$. This implies that for all $x \in \R\Lambda$, 
\begin{equation}
|f(z_{0}+x)|\leqslant C_{z_{0}} |e^{a(x,z_{0})}|=C_{z_{0}}e^{\Re a(x,z_{0})}.
\end{equation}

Let $\C_{\Lambda}^{m}$ be the maximal complex subspace of $\R\Lambda$. We need that the function $x\mapsto \Re (a(x,z_{0}))$ is bounded. In order to satisfy this, we will slightly modify the map $a$. There exists a unique linear map $g_{z_{0}}: \C_{\Lambda}^{m}\rightarrow \C$ such that $\Re(g_{z_{0}}(\lambda_{1}))  = \Re(a(\lambda_{1},z_{0}))$ for all $\lambda \in \Lambda$.

We then consider the summand of automorphy $\widetilde{a}_{z_{0}}: \Lambda \times \widetilde{U} \rightarrow \C$ defined by:
\[\widetilde{a}_{z_{0}}(\lambda, z) = a(\lambda, z) - g_{z_{0}}(z+\lambda)+g_{z_{0}}(z) = a(\lambda,z)-g_{z_{0}}(\lambda),\]
it defines the same line bundle as $a$. Hence, up to exchanging $a$ with $\widetilde{a}_{z_{0}}$, one can now assume that $\Re( a(\lambda_{1}, z_{0})) = 0$ for all $\lambda \in \Lambda$. For $x \in \C^{m}_{\Lambda}$ and $\lambda \in P$, one has:
\begin{eqnarray*}
a(x+\lambda_{1},z_{0}) &=& a(x,z_{0}+\lambda_{1})+a(\lambda_{1},z_{0})\\
& = &a(x,z_{0})+a(\lambda_{1},z_{0})
\end{eqnarray*}
(since by assumption on $a$, in its second argument it only depends on the last $n-m$  coordinates), so the restriction of $\Re a$ to $\C_{\Lambda}^{m}\times \{z_{0}\}$ is $\widetilde{T}.\Z^{2m}$-periodic therefore bounded. Consequently, the restriction of $f(z_{0}+ ~\cdot~)$ to $\C_{\Lambda}^{m}$ is a bounded holomorphic function, so it is constant and this holds for all $z_{0}\in\widetilde{U}$. This lets us see that $f$ does not depend on the first $m$ coordinates; we can see $f$ as a holomorphic function over $\widetilde{U}_{2}=\{z_{2}~|~z\in \widetilde{U}\}$. The zeros of $f$ are $(I_{n-m}~R)$-invariant, because of condition (\ref{condSurF}). Since $\C^{n}/\Lambda$ is a Cousin group, there is no $\sigma \in \Z^{n-m}\setminus\{0\}$ with $^{t}\sigma R\in\Z^{2m}$ so the group generated by $(I_{n-m}~R)$ is dense in $\R^{n-m}$. It is a consequence that if $f(w)=0$ for some $w\in \widetilde{U}_{2}$, then $f$ vanishes on $w+\R^{n-m}$ and also on the intersection of $\widetilde{U}_{2}$ and the complex linear subspace $\C^{n-m}$ generated by $w+\R^{n-m}$, i.e. $\widetilde{U}_{2}\cap\C^{n-m}=\widetilde{U}_{2}$. % tout entier par le principe des zéros isolés ($f$ est nulle sur un sous-ensemble de $\widetilde{U}_{2}$ dont la dimension réelle est supérieure ou égale à $\dim_{\C}\widetilde{U}_{2}$ donc nulle sur $\widetilde{U}_{2}$ tout entier).
Finally, if $s$ is a non-trivial section of $L$, it never vanishes and $L$ is holomorphically trivial.
}

Now we give a reformulation of the previous proposition which will be used in the next section:
\Lemme{\label{lemmeVogtPourOT}Let $U$ be an open set of a Cousin group $C \cong \C^{n}/\Lambda$ whose inverse image $\widetilde{U}$ in $\C^{n}$ is a convex domain and let $D$ be a positive divisor of $U$. Let $L$ be the line bundle over $U$ associated to $D$. Then $L$ is not topologically trivial, i.e. its first Chern class $c_{1}(L)\in H^{2}(U,\Z)$ is not zero.}
\Preuve{Assume that $L$ is topologically trivial. By theorem~\ref{propVogtGeneralisee}, since it admits a (non trivial) section $s_{1}$ given by $D$, it is holomorphically trivial. This means that it has a holomorphic section $s_{2}$ which never vanishes. The quotient $s_{1}/s_{2}$ is a non-constant holomorphic function on $U$ (since $s_{1}$ vanishes on $D$ and not on the complement of $D$), a contradiction. One can also notice that $s_{1}$ vanishes on $D$ so it is identically zero by the proof of the previous theorem, which leads to a contradiction too.}

\section{Hypersurfaces on OT-manifolds}
\noindent The goal of this section is to prove that the algebraic dimension of an OT-manifold is zero. We first briefly recall the construction of OT-manifolds.% and we prove this statement by using the results of the previous section, theorem~\ref{propVogtGeneralisee} in particular.

\subsection{Construction of OT-manifolds}
\noindent Let $K$ be a number field of degree $n$ over $\Q$; call $s$ (resp. $2t$) the number of real (resp. complex) embeddings of $K$, so that $n=s+2t$. For the construction, one needs to assume that the integers $s$ and $t$ are non-zero.

Then $\sigma(\mathcal{O}_{K})$ is a lattice of rank $n$ in $\C^{m}$ (where we set $m:=s+t$) hence we have a properly discontinuous action of $\sigma(\mathcal{O}_{K})$ on $\C^{m}$. 
%Considérons l'action multiplicative de $K$ sur $\C^{m}$, donnée par $$az:=(\sigma_{1}(a)z_{1}, ..., \sigma_{m}(a)z_{m})$$ pour $a\in K$ et $z\in \C^{m}$ (on remarque que $\sigma(\mathcal{O}_{K})$ est stable pour cette action). 
We also define $l: \mathcal{O}_{K}^{*,+} \rightarrow \R^{m}$ by 
$$l(a) = (\log |\sigma_{1}(a)|, ..., \log|\sigma_{s}(a)|, 2\log |\sigma_{s+1}(a)|, ..., 2\log|\sigma_{m}(a)|).$$ 
By Dirichlet's units theorem, $l(\mathcal{O}_{K}^{*,+})$ is a lattice in the vector space $L:= \{x \in \R^{m}~|~\sum_{i=1}^{m}x_{i}=0\}$. The projection $pr: L \rightarrow \R^{s}$ given by the $s$ first coordinates is surjective so there are subgroups $A$ of rank $s$ of $\mathcal{O}_{K}^{*,+}$ such that $pr(l(A))$ is a lattice of rank $s$ of $\R^{s}$. Such an $A$ is called \textbf{admissible}. We then look at the quotient $X:=X(K,A)=(\H^{s}\times\C^{t})/(A\ltimes \mathcal{O}_{K})$; it is a complex compact manifold of dimension $m$. This manifold is called an \textbf{Oeljeklaus-Toma manifold}, or \textbf{OT-manifold}.\\

We know (\cite{Oeljeklaus:2005aa}, lemma~2.4) that $\C^{m}/\sigma(\mathcal{O}_{K})$ admits no non-constant holomorphic function. In other words, the complex Lie group $C:=\C^{m}/\sigma(\mathcal{O}_{K})$ is a Cousin group.\\

We recall two lemmas and a definition (see \cite{Oeljeklaus:2005aa}) for further use:
\Lemme{Let $A$ be a subgroup of $\mathcal{O}_{K}^{*,+}$ which is not contained in $\Z$. Then the two following conditions are equivalent:
\begin{itemize}
\item[1.] The action of $A$ on $\mathcal{O}_{K}$ admits a (non-trivial) proper invariant sub-module of lower rank.
\item[2.] There exists a proper intermediate field extension $\Q \subset K' \subset K$ with $A \subset \mathcal{O}_{K'}^{*,+}$.
\end{itemize}}

%This lemma leads to the following definition:
\Definition{We say that $X(K,A)$ is of \textbf{simple type} if $A$ does not satisfy one of the equivalent conditions of the previous lemma.}

%Finally, we state a lemma which will be useful in the next section:
\Lemme{\label{lemme16OT}Let $\Q \subset K' \subset K$ an intermediate field extension with $A \subset \mathcal{O}_{K'}^{*,+}$ an admissible subgroup for $K$. Let $s', 2t'$ be the number of real and complex embeddings of $K'$ respectively. Then $s=s'$, $t'>0$ and $A$ is admissible for $K'$.}
 
For more details, see \cite{Oeljeklaus:2005aa}.

\subsection{Algebraic dimension}
\noindent We shall first consider the case of an OT-manifold of simple type. %Then we will use this result for the general case.

\Proposition{\label{dimAlgCasSimple}Let $X$ be an OT-manifold of simple type. Then it admits no irreducible divisor; in particular, the algebraic dimension of $X$ is zero.}
\Preuve{%We will show that $X$ admits no irreducible divisor. 
Assume that $D$ is an irreducible divisor of $X$ of multiplicity $1$.
We consider $X=X(K,A)$ as the quotient of $U:=(\H^{s}\times \C^{t})/\sigma(\mathcal{O}_{K})$ (which is diffeomorphic to $(\R_{>0})^{s}\times (\S^{1})^{n}$) by $A\cong \Z^{s}$.
One has the following commutative diagram: 
\begin{equation}
\label{diagrammeOTCousin1}
\xymatrix{
U \ar[rrr]^{\displaystyle A \cong \Z^{s}}_{p} \ar[dd]_{\displaystyle  (\S^{1})^{s+2t}}^{q} & & & X \ar[dd]^{\displaystyle (\S^{1})^{s+2t}}_{q'}\\
& & & \\
(\R_{>0})^{s} 
%\ar[dd]_{\displaystyle \Z}
 \ar[rrr]_{\displaystyle \Z^{s}}^{p'}& & & (\S^{1})^{s}
%\ar[dd]^{\displaystyle \bar\Z}
%& & & \\
%Y & & & W 
}\end{equation}
Because $X$ is of simple type, one has $\dim H^{2}(X)=\displaystyle \binom{s}{2}=\dim H^{2}((\S^{1})^{s})$ by proposition~2.3 of \cite{Oeljeklaus:2005aa}. Thus the first Chern class $c_{1}(L(D))$ of the line bundle $L(D)$ over $X$ associated to $D$ is the pull-back by $q'$ of an element $\omega \in H^{2}((\S^{1})^{s},\Z)$. The fact that the map $(q')^{*}: H^{2}((\S^{1})^{s},\Z)\rightarrow H^{2}(X,\Z)$ is injective is established in the proof of proposition~2.3 of \cite{Oeljeklaus:2005aa}. The commutativity of diagram~(\ref{diagrammeOTCousin1}) implies that the pull-back $\widetilde{\omega}$ over $U$ of $\omega$ by $q' \circ p$ is the same as the pull-back by $p' \circ q$ of $\omega$. Furthermore, one has $(p')^{*}(\omega)=0$ since $H^{2}((\R_{>0})^{s},\Z)=0$. Hence $\widetilde{\omega}=0$. This is a contradiction to lemma~\ref{lemmeVogtPourOT}, since $\widetilde{\omega}=c_{1}(\widetilde{L}(D))$ where $\widetilde{L}(D)$ is the bundle associated to the divisor $p^{-1}(D)$ of $U$. 
}

Now we can omit the hypothesis that $X$ is of simple type:

\Theoreme{\label{dimAlgCasGeneral}Let $X$ be an OT-manifold. Then there are no complex-analytic hypersurfaces on $X$, in particular the algebraic dimension of $X$ is zero.}
\Preuve{By proposition~\ref{dimAlgCasSimple}, we can assume that $X$ is not of simple type. Consequently there exists a field extension $\Q \subset K' \subset K$ with $A \subset \mathcal{O}_{K'}^{*,+}$ and we call $s'$ (respectively, $2t'$) the number of real (respectively, complex) embeddings of $K'$. Since $A$ is admissible for $K$, we have $s=s'$ and $A$ is admissible for $K'$ (see lemma~\ref{lemme16OT}). Without loss of generality, we can suppose that $K'$ is the smallest subfield of $K$ such that $A \subset \mathcal{O}_{K'}^{*,+}$, i.e. $X(K',A)$ is of simple type. Call $\sigma_{1}, ..., \sigma_{s}, \sigma_{s+1}, ..., \sigma_{s+t},\bar\sigma_{s+1}, ..., \bar\sigma_{s+t}$ the $s+2t$ embeddings of $K$ and $\sigma'_{1}, ..., \sigma'_{s}, \sigma'_{s+1}, ..., \sigma'_{s+t'},\bar\sigma'_{s+1}, ..., \bar\sigma'_{s+t'}$ the $s+2t'$ embeddings of $K'$.\\

As before, we note $$\application{\sigma}{\mathcal{O}_{K}}{\C^{s+t}}{a}{(\sigma_{1}(a), ..., \sigma_{s+t}(a)).}$$% et $$\application{\sigma'}{\mathcal{O}_{K'}}{\C^{s+t'}}{a}{(\sigma'_{1}(a), ..., \sigma'_{s+t'}(a)).}$$

We look at the complex linear subspace $V_{K'}$ of $\C^{s+t}$ spanned by $\sigma(\mathcal{O}_{K'})$. We denote by $C$ (respectively, $C'$) the Cousin group $\C^{s+t} / \sigma(\mathcal{O}_{K})$ (respectively, $V_{K'} / \sigma(\mathcal{O}_{K'})$). The group $C'$ is a closed complex Lie subgroup of $C$. We study the two open sets $U:=(\H^{s}\times \C^{t})/\sigma(\mathcal{O}_{K})\subset C$ and $U':=((\H^{s}\times \C^{t})\cap V_{K'})/\sigma(\mathcal{O}_{K'})\subset C'$. The quotients $U/A$ and $U'/A$ are OT-manifolds that we call $X$ and $X'$ respectively; moreover, $X'$ is a compact complex submanifold of $X$.\\

We consider the following commutative diagram:
\begin{equation}
\label{diagrammeCommutatifOTCousin}
\xymatrix{
0 \ar[r] & C'~ \ar@{^{`}->}[r] & C \ar@{->>}[r]^{q~~~}  &  C/C' \ar[r] & 0\\
& U' \vphantom{\int}~ \ar@{^{`}->}[r]\ar@{^{`}->}[u] \ar[d]_{p'} & U \vphantom{\int} \ar@{^{`}->}[u]\ar[d]^{p\phantom{'}}\ar[r] & q(U)\vphantom{\int}\ar@{^{`}->}[u] & \\
& X'~ \ar@{^{`}->}[r]^{i} & X & &
}
\end{equation}

Let $D$ be a irreducible divisor of multiplicity $1$ of $X$ and let $L(D)$ be the associated holomorphic line bundle. By restriction we get a bundle $L':= i^{*}L(D)$ over $X'$. 

Since $X'$ is of simple type, the line bundle $(p')^{*}(L')$ over $U'$ is topologically trivial, see the proof of proposition~\ref{dimAlgCasSimple}.
The maximal complex subgroup $H'\cong \C^{t'}$ of the maximal compact torus $T' \cong (\S^{1})^{s+2t'}$ of $C'$ is a subgroup of the maximal complex subgroup $H \cong \C^{t}$ of the maximal compact torus $T\cong (\S^{1})^{s+2t}$ of $C$. We denote by $B\cong \C^{t-t'}$ a connected complex subgroup of $T$ such that $B \times H' \cong H$. Since $B$ is a subgroup of $T$ it acts on $U$ and this action is transitive on the leaves of the $U'$-foliation of $U$ induced by the $C'$-foliation of $C$. The induced action of $B$ on $H^{2}(U,\Z)$ is trivial since $H^{2}(U,\Z)$ is discrete and $B$ is connected, i.e. we have $b^{*}c_{1}(L) = c_{1}(L)$ for every line bundle $L$ above $U$ and all $b\in B$.

The inverse image $\widetilde{D}:=p^{-1}(D)$ is a divisor of $U$. We shall now prove that the divisor $\widetilde{D}$ is saturated by the leaves of the $U'$-foliation of $U$.

\noindent We have 
\begin{equation}
\label{etoile}c_{1}(L(\widetilde{D}))=p^{*}(c_{1}(L(D)))\in H^{2}(U,\Z).
\end{equation}
Since the diagram~(\ref{diagrammeCommutatifOTCousin}) commutes, we have $c_{1}(L(\widetilde{D})|_{U'})=c_{1}(L(\widetilde{D}))|_{U'}=c_{1}(p^{*}L')=0=c_{1}(b^{*}L(\widetilde{D}))|_{U'}$. There are three possible cases for the intersection $b(\widetilde{D}) \cap U'$. Either this intersection is $U'$, or it is empty (these two cases fit our purpose), or it is a divisor of $U'$. The last case can not occur since $c_{1}(b^{*}L(\widetilde{D}))|_{U'}$ would be non-zero by lemma~\ref{lemmeVogtPourOT}, which is a contradiction to (\ref{etoile}). Hence $\widetilde{D}$ is saturated by the leaves of the $U'$-foliation of $U$.%~\\~\\truc avec le groupe B pour saturer D par U'~\\~\\

Now, $\widetilde{D}$ induces a divisor $q(\widetilde{D})$ of $q(U)$. One has $q(U)=C/C'$ because the leaves of the $U'$-foliation of $U$ are in bijection with those of the $C'$-foliation of $C$. Hence, $q(\widetilde{D})$ lifts as a divisor in $C$, which is stable under the action of the group $A\ltimes C'$. We still denote by $\widetilde{D}$ this divisor. 

Since the group $\mathcal{O}_{K}^{*,+}$\! is abelian, the action of $\mathcal{O}_{K}^{*,+}$\! on $U$ induces an action of $\mathcal{O}_{K}^{*,+}$\! on $X$. Let $\eta \in \mathcal{O}_{K}^{*,+}$\!, we have $\widetilde{\eta D} = \eta \widetilde{D}$. We use the $C'$-invariance of $\widetilde{\eta D}$ to write $\eta\widetilde{D}+C'=\eta\widetilde{D}$ hence $\widetilde{D}$ is both $C'$-invariant and $\eta^{-1}C'$-invariant. 
%\textcolor{red}{Let $f\in\mathcal{M}(X)$ be a meromorphic function on $X$, then $f$ lifts as a meromorphic function on $U$ invariant by $A$ (still denoted by $f \in \mathcal{M}^{A}(U)$). The ``intermediate result'' of the previous paragraph tells us that the level sets of $f$ are saturated by the $U'$-foliation of $U$ so $f$ is constant on the leafs of the $U'$-foliation of $U$. }
%\textbf{truc à voir: } $U' \hookrightarrow U \rightarrow C/C' $ la dernière flèche est surjective 
%Therefore, $f$ induces an application on $q(U) = C/C'$ (this equality comes from the fact that the leaves of the $U'$-foliation of $U$ are in bijection with those of the $C'$-foliation of $C$) so it lifts as a meromorphic function on $C$, which is invariant by the group $A\ltimes C'$; we still denote by $f \in \mathcal{M}^{A\ltimes C'}(C)$ this function; we have in fact shown that $\mathcal{M}(X)=\mathcal{M}^{A\ltimes C'}(C)$. 

%Since the group $\mathcal{O}_{K}^{*,+}$ is abelian, the action of $\mathcal{O}_{K}^{*,+}$ on $U$ induces an action of $\mathcal{O}_{K}^{*,+}$ on $X$; the action of $\mathcal{O}_{K}^{*,+}$ on $\mathcal{M}(X)$ corresponds to the one on $\mathcal{M}^{A\ltimes C'}(C)$. 
%We see $f$ as an element of $\mathcal{M}^{A\ltimes C'}(C)$: for all $\eta\in \mathcal{O}_{K}^{*,+}$ we set $f_{\eta}(x):= f(\eta^{-1}x)$. Then for all $x \in C$, we have $f_{\eta}(x+C') = f_{\eta}(x)$ i.e. $f$ is $C'$-invariant and $\eta^{-1}C'$-invariant. 
Let $J$ be the smallest connected complex subgroup of $C$ containing $C'$ and $\eta C'$ for all $\eta\in \mathcal{O}_{K}^{*,+}$. If $J$ is a proper subgroup of $C$, its maximal compact subgroup is defined by a sublattice of $\mathcal{O}_{K}$ which is stable under the action of $\mathcal{O}_{K}^{*,+}$. This is impossible because there is no sublattice of $\mathcal{O}_{K}$ stable by $\mathcal{O}_{K}^{*,+}$. Indeed, if it were the case, since $\mathcal{O}_{K}^{*,+}$ contains a primitive element of $K$, this element should have a minimal polynomial with degree strictly smaller than the one of $K$ and this is impossible. Finally, $J=C$ and $\widetilde{D}$ is invariant by $C$, which is a contradiction. Hence, $X$ admits no divisor and its algebraic dimension is zero.}

\section{A special class of Cousin groups}
\subsection{Introduction}
\noindent In \cite{Vogt:1982aa}, Vogt exhibits a special class of Cousin groups by prooving\footnote{ The theorem given in \cite{Vogt:1982aa}, p.~208 has $8$ equivalent assertions, here we only recall three of them.}:
\Theoreme{\label{thCaractClasseCousinVogt}Let $C=\C^{n}/\Lambda$ be a Cousin group. Then the following conditions are equivalent:
\begin{itemize}
%\item[1.\phantom{'}] Chaque factor of automorphy $\alpha: \Lambda \times \C^{n}\rightarrow \C$ est équivalent à un facteur theta\footnote{Un factor of automorphy $\alpha$ est un \emph{facteur theta} si pour chaque $\lambda\in\Lambda$ il existe un polynôme de degré $1$ noté $a_{\lambda}$ such that $\alpha_{\lambda}=\exp(a_{\lambda})$.}.
%\item[2.\phantom{'}] Tout summand of automorphy $a: \Lambda \times \C^{n} \rightarrow \C$ est équivalent à un homomorphisme $\Lambda \rightarrow \C$.
%\item[2'.] Tout fibré en droite topologiquement trivial sur $X$ est donné par une représentation de $\Lambda$.
%\item[3.\phantom{'}] Tout fibré en droite topologiquement trivial sur $X$ est homogène.
%\item[4.\phantom{'}] Pour toute $(0,1)$-forme $\omega$ sur $\C^{n}$ qui est $\Lambda$-périodique et $\bar\partial$-fermée, il existe un fonction $f$ différentiable sur $\C^{n}$, $\Lambda$-périodique et des constantes $a_{1}, ..., a_{n}\in\C$ telles que $\omega = \bar\partial f + \sum_{j=1}^{n}a_{j}d\bar z_{j}$.
\item[1.] The space $H^{1}(C,\mathcal{O})$ is finite-dimensional.
%\item[6.\phantom{'}] L'application naturelle $H^{1}(X,\C)\rightarrow H^{1}(X,\mathcal{O})$ est surjective.
%\item[7.\phantom{'}] Soit $P$ une base de périodes pour $\Lambda$ et $T:=\C^{m}$...
\item[2.] Let $P=(I_{n}~S)$ be a period basis of $\Lambda$. Then there exist constants $\mathcal{C}>0$ and $a\geqslant 0$ such that $\|^{t}\sigma S+\,^{t}\tau\|\geqslant \mathcal{C} \exp(-a|\sigma|)$ for all $\sigma \in \Z^{n}\setminus\{0\}$ and all $\tau\in\Z^{m}$, where $n+m$ is the rank of $\Lambda$.
%\item[9.\phantom{'}] Soit $P=\left(\begin{array}{lr} 0& T^{*}\\I_{n-m}& R\end{array}\right)$ une base de périodes pour $\Lambda$ vérifiant les conditions de la proposition~\ref{formesNormalesReseauCousin} et $T:=\C^{m}/T^{*}\Z^{2m}$ le tore de dimension $m$ correspondant à $T^{*}$. Alors l'application $H^{1}(T,\mathcal{O})\rightarrow H^{1}(X,\mathcal{O})$ induite par $p: X \rightarrow T$ est surjective.
\item[3.] Every line bundle over $C$ comes from a theta factor.
\end{itemize}}

\noindent Let $\alpha$ be a real irrational algebraic number. In his paper, Vogt gives the following example of a lattice in $\C^{2}$ given by $P=\left(\begin{array}{ccc}0& 1 & i\\1 & \alpha & 0\end{array}\right)$ which defines a Cousin group satisfying condition n°2. Here we give a general theorem in this setting which in particular applies to all Cousin groups appearing in the construction of OT-manifolds.%The goal of this section is to show that if a Cousin group is defined by a lattice whose vectors have algebraic numbers as coordinates, it belongs to this class (because it satisfies condition n°2). As a consequence of this result, the Cousin groups appearing in the construction of OT-manifolds belong to this class of Cousin groups. 

%\subsection{Preliminaries}
\noindent We shall use the following generalization of Liouville's theorem which can be found in \cite{Feldman:1998aa} (theorem~1.5, page~27):

\Theoreme{\label{thLiouvGen}Let $\alpha_{1}, ..., \alpha_{m}$ be algebraic numbers, of respective degrees $n_{k}$, with $\operatorname{deg} \Q(\alpha_{1}, ..., \alpha_{m})=n$, and let \[P(z_{1}, ..., z_{m})=\sum_{k_{1}=0}^{N_{1}}\cdots\sum_{k_{m}=0}^{N_{m}}a_{k_{1}, ..., k_{m}}z_{1}^{k_{1}}\cdots z_{m}^{k_{m}}\in\Z[z_{1}, ..., z_{m}].\]
If $P(\alpha_{1}, ..., \alpha_{m})$ is non-zero, then one has the inequality
\[|P(\alpha_{1}, ..., \alpha_{m})| \geqslant L(P)^{1-\delta n}\prod_{k=1}^{m}L(\alpha_{k})^{-\delta N_{k}n/n_{k}},\]
with $\delta=1$ if all the $\alpha_{i}$ are real, $\delta=1/2$ otherwise and where $L(P)$ is the sum of the absolute values of the coefficients of $P$ (and $L(\alpha)$ is the quantity $L(\mu)$, with $\mu$ being the minimal polynomial of $\alpha$).}

\subsection{A class of Cousin groups satisfying condition n°2}
We shall show:
\Theoreme{\label{grpeCousinAlgebriqueCondition2}Let $\Lambda \subset \C^{n}$ be a lattice such that $C=\C^{n}/\Lambda$ is a Cousin group with a period basis whose coefficients are all algebraic numbers, then $C$ satisfies the equivalent conditions of theorem~\ref{thCaractClasseCousinVogt}.}
\Preuve{%Choose a basis $(x_{1}, ..., x_{n+m})$ of $\Lambda$ whose coefficients are algebraic numbers, i.e. %, which we write as vectors: 
%
%\[\left(\begin{array}{c}\alpha_{1,1} \\ \vdots \\\alpha_{1,n}\end{array}\right), ..., \left(\begin{array}{c}\alpha_{n+m,1} \\ \vdots \\\alpha_{n+m,n}\end{array}\right).\]
%
%\textcolor{red}{directement p = i s}
%Suppose $\Lambda$ is given by a period basis of $X$:
%\[P=\left(\begin{array}{c}\alpha_{1,1} \\ \vdots \\\alpha_{1,n}\end{array} ... \begin{array}{c}\alpha_{n+m,1} \\ \vdots \\\alpha_{n+m,n}\end{array}\right),\]
%where the $\alpha$
%Let $A$ be the matrix defined as follows: 
%\[A:=\left(\begin{array}{c}\alpha_{1,1} \\ \vdots \\\alpha_{1,n}\end{array} ... \begin{array}{c}\alpha_{n,1} \\ \vdots \\\alpha_{n,n}\end{array}\right).\]
%
%This matrix is invertible and the matrix $A^{-1}P$ defines a Cousin group isomorphic to $X$. The new period basis is then 
%\[P':=A^{-1}P = (I_{n}~S)\]
Let us write a period basis of $\Lambda$ as a matrix \[(I_{n}~S)\]
where $I_{n}$ is the identity matrix of size $n$ and $S=(s_{i,j})$ is a matrix with $n$ rows and $m$ columns whose coefficients are algebraic numbers over $\Q$.

In order the check that the second condition of theorem~\ref{thCaractClasseCousinVogt} is satisfied, we have to verify that there exist constants $\mathcal{C}>0$ and $a\geqslant 0$ such that for all $\mu \in \Z^{n}\setminus\{0\}$ and all $\nu \in \Z^{m}$, one has
\[\|^{t}\mu S +\,^{t}\nu\|\geqslant \mathcal{C}\exp(-a|\mu|).\]

It is enough to prove this inequality for any non-zero coordinate of the vector $^{t}\mu S +\,^{t}\nu$ and then consider the infinity norm of this vector. %, which is non-zero because we know that $^{t}\mu S \not \in \Z^{m}$. 
Denote $^{t}\mu=(\mu_{1}, ..., \mu_{n})$ and $^{t}\nu=(\nu_{1}, ..., \nu_{m})$. Then the $k$-th coordinate of the vector $^{t}\mu S +\,^{t}\nu$ (for $k \in\{1, ..., m\}$) is 
\[\mu_{1}s_{1,k}+\cdots + \mu_{n}s_{n,k}+\nu_{k}.\]
We now distinguish the two following cases: either we have the inequality $|\nu_{k}|\leqslant 2 |\mu_{1}s_{1,k}+\cdots+ \mu_{n}s_{n,k}|$, or its converse holds.\\

Suppose first $|\nu_{k}|\leqslant 2 |\mu_{1}s_{1,k}+\cdots+ \mu_{n}s_{n,k}|$. By theorem~\ref{thLiouvGen}, we have 
\[|\mu_{1}s_{1,k}+\cdots + \mu_{n}s_{n,k}+\nu_{k}| \geqslant \mathcal{C} (|\mu_{1}|+\cdots+|\mu_{n}|+|\nu_{k}|)^{p},\]
where $\mathcal{C}>0$ and $p<0$ are constants independent of $\mu$ and $\nu$. Recall here that we assumed that $\mu_{1}s_{1,k}+\cdots + \mu_{n}s_{n,k}+\nu_{k}$ is non-zero.
By assumption on $\nu_{k}$, we have 
\begin{eqnarray*}
(|\mu_{1}|+\cdots+|\mu_{n}|+|\nu_{k}|)^{p} &\geqslant& (|\mu_{1}|+\cdots+|\mu_{n}|+2 |\mu_{1}s_{1,k}+\cdots+ \mu_{n}s_{n,k}|)^{p}\\
&\geqslant & \mathcal{C}'(|\mu_{1}|+\cdots +|\mu_{n}|)^{p}\\
&\geqslant& \mathcal{C}'\exp(-|p| (|\mu_{1}|+\cdots +|\mu_{n}|)),
\end{eqnarray*}
with $\mathcal{C}'$ independent of $\mu$ and $\nu$. \\%This is what we wanted to prove.\\

Now suppose that $|\nu_{k}|> 2 |\mu_{1}s_{1,k}+\cdots+ \mu_{n}s_{n,k}|$. We write the reverse triangle inequality: 
\begin{eqnarray*}|\mu_{1}s_{1,k}+\cdots + \mu_{n}s_{n,k}+\nu_{k}| &\geqslant& \left|\vphantom{\sum}|\nu_{k}|-|\mu_{1}s_{1,k}+\cdots + \mu_{n}s_{n,k}|\right|\\
&=&|\nu_{k}|-|\mu_{1}s_{1,k}+\cdots + \mu_{n}s_{n,k}|\\
&\geqslant& |\mu_{1}s_{1,k}+\cdots + \mu_{n}s_{n,k}|.
\end{eqnarray*}
If $\mu_{1}s_{1,k}+\cdots + \mu_{n}s_{n,k}$ vanishes, we have $|\nu_{k}|\geqslant 1$ and the result is obtained; otherwise, we can again use theorem~\ref{thLiouvGen}: 
\begin{eqnarray*}
|\mu_{1}s_{1,k}+\cdots + \mu_{n}s_{n,k}|&\geqslant&\mathcal{C}''(|\mu_{1}|+\cdots+|\mu_{n}|)^{q}\\
&\geqslant&\mathcal{C}''\exp(-|q|(|\mu_{1}|+\cdots|\mu_{n}|))
\end{eqnarray*}
where $\mathcal{C}''>0$ and $q<0$ are constants which do not depend on $\mu$ and $\nu$. 
}

%La preuve ci-dessus permet en fait d'établir le résultat (plus fort) suivant: 
%\Theoreme{Soit $\Lambda \subset \C^{n}$ un réseau such that $X=\C^{n}/\Lambda$ soit un Cousin group et admettant une base de périodes dont tous les coefficients sont des nombres algébriques, alors $X$ vérifie la condition $2$ du theorem~\ref{thCaractClasseCousinVogt}.}

\noindent \textbf{Application}\\
\noindent We recall that in the construction of an OT-manifold, the group $C:=\C^{m}/\sigma(\mathcal{O}_{K})$ is a Cousin group. As a corollary of theorem~\ref{grpeCousinAlgebriqueCondition2}, we have: 
\Corollaire{The Cousin group $C$ satisfies condition n°2 of theorem~\ref{thCaractClasseCousinVogt}.}

\bibliographystyle{amsplain}
\bibliography{bibliographie}

\end{document}